# New generalization of contraction mapping by new control function

S. M. A. Aleomraninejad

Department of Mathematics, Qom University of Technology, Qom, Iran

email: aleomran63@yahoo.com

**Abstract.** In this paper, we introduce the new generalization of contraction mapping by a new control function and an altering distance . We establish some existence results of fixed point for such mappings. Our results reproduce several old and new results in the literature.

## 1 Introduction

The first important result obtained on fixed points for contractive-type mappings was the well-known Banach contraction theorem, published for the first time in 1922 ([2]). In the general setting of complete metric spaces, this theorem runs as follows.

**Theorem 1.1.** *Let $(X, d)$ be a complete metric space, $\beta \in (0,1)$ and let $T : X \to X$ be a mapping such that for each $x, y \in X$,*

$$d(Tx, Ty) \leq \beta d(x, y).$$

*Then $T$ has a unique fixed point $a \in X$ such that for each $x \in X$, $\lim_{n \to \infty} T^n x = a$.*

In order to generalize this theorem, several authors have introduced various types of contraction inequalities. In 2002 Branciari proved the following result (see [4]).

**Theorem 1.2.** *Let $(X, d)$ be a complete metric space, $\beta \in (0,1)$ and $T : X \longrightarrow X$ a mapping such that for each $x, y \in X$,*

$$\int_0^{d(Tx,Ty)} f(t) dt \leq \beta \int_0^{d(x,y)} f(t) dt,$$

*where $f : [0, \infty) \to (0, \infty)$ is a Lebesgue integrable mapping which is summable (i.e., with finite integral on each compact subset of $[0, \infty)$) and for each $\varepsilon > 0$, $\int_0^{\varepsilon} f(t) dt > 0$. Then $T$ has a unique fixed point $a \in X$ such that for each $x \in X$, $\lim_{n \to \infty} T^n x = a$.*

Rhoades [13] and Djoudi et al. [6] extended the result of Branciari and proved the following fixed point theorems.

**Theorem 1.3.** *[13] Let $(X,d)$ be a complete metric space, $k \in [0,1)$, $T : X \to X$ a mapping satisfying for each $x, y \in X$,*

$$\int_0^{d(Tx,Ty)} \varphi(t)dt \leq k \int_0^{M(x,y)} \varphi(t)dt$$

*where*

$$M(x,y) = \max\{d(x,y), d(x,Tx), d(y,Ty), \frac{1}{2}[d(x,Ty) + d(y,Tx)]\}$$

*and $\varphi : \mathbb{R}^+ \to \mathbb{R}^+$ be as in theorem 1.2. Then $T$ has a unique fixed point $x \in X$.*

**Theorem 1.4.** *[6] Let $(X,d)$ be a complete metric space and $T : X \to X$ a mapping satisfying for each $x, y \in X$,*

$$\int_0^{d(Tx,Ty)} \varphi(t)dt \leq h(\int_0^{M(x,y)} \varphi(t)dt)$$

*where*
$$M(x,y) = \max\{d(x,y), d(x,Tx), d(y,Ty), d(x,Ty), d(y,Tx)\},$$

$h : \mathbb{R}^+ \to \mathbb{R}^+$ *is subadditive, nondecreasing and continuous from the right such that $h(t) < t$, for all $t > 0$ and $\varphi : \mathbb{R}^+ \to \mathbb{R}^+$ be as in theorem 1.2. Then $T$ has a unique fixed point $x \in X$.*

In 1984, M.S. Khan, M. Swalech and S. Sessa [9] expanded the research of the metric fixed point theory to the category $\Psi$ by introducing a new function which they called an altering distance function. For $\psi : \mathbb{R}^+ \to \mathbb{R}^+$ we say that $\psi \in \Psi$ if

1. $\psi(t) = 0$ *if and only if* $t = 0$,

2. $\psi$ is monotonically non-decreasing,

3. $\psi$ is continuous.

The following lemma shows that contractive conditions of integral type can be interpreted as contractive conditions involving an altering distance.

**Lemma 1.5.** *Let $\varphi : \mathbb{R}^+ \to \mathbb{R}^+$ be as in Theorem 1.2. Define $\psi(t) = \int_0^t \varphi(\tau)d\tau$, for $t \in \mathbb{R}^+$. Then $\psi$ is an altering distance.*

Khan et al. using this altering distance to extend the Banach Contraction Principle as follows:

**Theorem 1.6.** *[9] Let $(X,d)$ be a complete metric space, $\beta \in (0,1)$ and $T : X \longrightarrow X$ a mapping such that for each $x, y \in X$,*

$$\psi[d(Tx,Ty)] \leq \beta \psi[d(x,y)]$$

*where $\psi \in \Psi$. Then $T$ has a unique fixed point $a \in X$ such that for each $x \in X$, $\lim_{n \to \infty} T^n x = a$.*

It is easy to see that if $\psi(t) = t$, we obtain the Banach Contraction Principle and by lemma 1.5, we obtain theorem 1.2. Dutta et al. [8], Dori [7], Choudhury et al. [5] and Morals et al. [11] extended the results of Khan and proved the following fixed point theorems.

**Theorem 1.7.** *[8] Let $(X, d)$ be a complete metric space and let $T : X \to X$ be a mapping satisfying*

$$\psi(d(Tx, Ty)) \leq \psi(d(x, y)) - h(d(x, y))$$

*for each $x, y \in X$, where $\psi, h : \mathbb{R}^+ \to \mathbb{R}^+$ are continuous and non-decreasing function such that $\psi(t) = h(t) = 0$ if and only if $t = 0$. Then $T$ has a unique fixed point $x \in X$.*

**Theorem 1.8.** *[7] Let $(X, d)$ be a complete metric space and let $T : X \to X$ be a mapping satisfying*

$$\psi(d(Tx, Ty)) \leq \psi(M(x, y)) - h(M(x, y)))$$

*for each $x, y \in X$, where*

$$M(x, y) = \max\{d(x, y), d(x, Tx), d(y, Ty), \frac{1}{2}[d(x, Ty) + d(y, Tx)]\},$$

*$\psi \in \Psi$ and $h : \mathbb{R}^+ \to \mathbb{R}^+$ is a lower semi-continuous function such that $h(t) = 0$ if and only if $t = 0$. Then $T$ has a unique fixed point $x \in X$.*

**Theorem 1.9.** *[5] Let $(X, d)$ be a complete metric space and let $T : X \to X$ be a mapping satisfying*

$$\psi(d(Tx, Ty)) \leq \psi(M(x, y)) - h(\max\{d(x, y), d(y, Ty)\})$$

*for each $x, y \in X$, where*

$$M(x, y) = \max\{d(x, y), d(x, Tx), d(y, Ty), \frac{1}{2}[d(x, Ty) + d(y, Tx)]\},$$

*$\psi \in \Psi$ and $h : \mathbb{R}^+ \to \mathbb{R}^+$ is a continuous function such that $h(t) = 0$ if and only if $t = 0$. Then $T$ has a unique fixed point $x \in X$.*

**Theorem 1.10.** *[11] Let $(X, d)$ be a complete metric space and $T : X \longrightarrow X$ a mapping which satisfies the following condition:*

$$\psi[d(Tx, Ty)] \leq a\psi[d(x, y)] + b\psi[m(x, y)]$$

*for all $x, y \in X$, $a > 0$, $b > 0$, $a + b < 1$ where*

$$m(x, y) = d(y, Ty)\frac{1 + d(x, Tx)}{1 + d(x, y)}$$

*for all $x, y \in X$. where $\psi \in \Psi$. Then $T$ has a unique fixed point $a \in X$ such that for each $x \in X$, $\lim_{n \to \infty} T^n x = a$.*

On the other hand, in 2008, Suzuki introduced a new method in [14] and then his method was extended by some authors (see for example [1], [10], [12]). The aim of this paper is to provide a new and more general condition for $T$ which guarantees the existence of its fixed point. Our results generalize several old and new results in the literature. In this way, consider $\Phi$ the set of all control function $\phi : [0, \infty)^k \longrightarrow [0, \infty)$ satisfying
(i) $\phi(0, 0, ..., 0) = 0$,
(ii) $\lim_{n \to \infty} \phi(t_{1n}, t_{2n}, ..., t_{kn}) \leq \phi(t_1, t_2, ..., t_k)$ whenever $(t_{1n}, t_{2n}, ..., t_{kn}) \to (t_1, t_2, ..., t_k)$,
and $R$ the set of all continuous function $g : [0, \infty)^5 \longrightarrow [0, \infty)$ satisfying the following conditions:
(i) $g(1, 1, 1, 0, 2), g(1, 1, 1, 1, 1) \in (0, 1]$,
(ii) $g$ is subhomogeneous, i.e.
$g(\alpha x_1, \alpha x_2, \alpha x_3, \alpha x_4, \alpha x_5) \leq \alpha g(x_1, x_2, x_3, x_4, x_5)$ for all $\alpha \geq 0$.
(iii) if $x_i, y_i \in [0, \infty)$, $x_i \leq y_i$ for $i = 1, ..., 5$ we have $g(x_1, x_2, x_3, x_4, x_5) \leq g(y_1, y_2, y_3, y_4, x_5)$

**Example 1.1.** *Define $g(x_1, x_2, x_3, x_4, x_5) = \frac{1}{2} \max\{x_i\}_{i=1}^5$. It is obvious that $g \in R$.*

**Example 1.2.** *Define $g(x_1, x_2, x_3, x_4, x_5) = \max\{x_1, x_2, x_3, \frac{x_4+x_5}{2}\}$. It is obvious that $g \in R$.*

**Proposition 1.11.** *If $g \in R$ and $u, v \in [0, \infty)$ are such that*
$$u < \max\{g(v, v, u, v, u), g(v, u, v, v + u, 0)\},$$
*then $u < v$.*

*Proof.* Without loss of generality, we can suppose $u < g(v, u, v, v + u, 0)$. If $v \leq u$, then
$$u < g(v, u, v, v + u, 0) \leq g(u, u, u, 2u, 0) \leq ug(1, 1, 1, 2, 0) \leq u$$
which is a contradiction. Thus $u < v$. $\square$

**Lemma 1.12.** *Let $\psi \in \Psi$ and $\phi \in \Phi$ such that for every $t_i \in \mathbb{R}^+$,*
$$\phi(t_1, t_2, ..., t_k) < \psi(\max_{i=1,...,k} t_i).$$
*If for $t, s_i \in \mathbb{R}^+$ we have*
$$\psi(t) \leq \phi(s_1, s_2, ..., s_k),$$
*then*
$$t < \max_{i=1,...,k} s_i.$$

*Proof.* Let $S = \max_{i=1,...,k} s_i$. Suppose that $t \geq S$. Then
$$\psi(t) \geq \psi(S) > \phi(s_1, s_2, ..., s_k),$$
which is a contradiction. $\square$

**Lemma 1.13.** *Suppose that $\{s_n\}$ be a sequence of non-negative real numbers such that $s_{n+1} \leq s_n$. Then $s_n$ is convergent.*

**Lemma 1.14.** *[2] Let $(X,d)$ be a metric space and $\{x_n\}$ be a sequence in $X$ such that*

$$\lim_{n \to \infty} d(x_n, x_{n+1}) = 0.$$

*If $\{x_n\}$ is not a Cauchy sequence in $X$, then there exist an $\varepsilon_0 > 0$ and sequences of positive integers $m_k$ and $n_k$ with $m_k > n_k > k$ such that*

$$d(x_{m_k}, x_{n_k}) \geq \varepsilon_0, \quad d(x_{m_k-1}, x_{n_k}) < \varepsilon_0$$

*and*

*(i)* $\lim_{k \to \infty} d(x_{m_k-1}, x_{n_k+1}) = \varepsilon_0$,
*(ii)* $\lim_{k \to \infty} d(x_{m_k}, x_{n_k}) = \varepsilon_0$,
*(iii)* $\lim_{k \to \infty} d(x_{m_k-1}, x_{n_k}) = \varepsilon_0$.

## 2 Main Results

The following theorem is the main result of this paper.

**Theorem 2.1.** *Let $(X,d)$ be a complete metric space, $T : X \longrightarrow X$ a mapping, $\alpha \in (0, \frac{1}{2}]$, $\psi \in \Psi$ and $\phi \in \Phi$ such that for every $t_i \in \mathbb{R}^+$ with $(t_1, t_2, ..., t_k) \neq (0, 0, ..., 0)$,*

$$\phi(t_1, t_2, .., t_k) < \psi(\max_{i=1,...,k} t_i).$$

*Suppose that $\{g_i\}_{i=1}^k$ be a sequence in $R$ and $\alpha d(x, Tx) \leq d(x, y)$ implies*

$$\psi[d(Tx, Ty)] \leq \phi(g_1(M_{xy}), g_2(M_{xy}), ..., g_k(M_{xy}))$$

*for all $x, y \in X$, where*

$$M_{xy} = (d(x, y), d(y, Ty), d(x, Tx), d(x, Ty), d(y, Tx))$$

*for all $x, y \in X$. Then $T$ has a unique fixed point in $X$.*

*Proof.* Fix arbitrary $x_0 \in X$ and let $x_1 = Tx_0$. We have $\alpha d(x_0, Tx_0) < d(x_0, x_1)$. Hence,
$$\psi[d(Tx_0, Tx_1)] \leq \phi(g_1(M_{x_0 x_1}), g_2(M_{x_0 x_1}), ..., g_k(M_{x_0 x_1})).$$

Then by lemma 1.12 we have

$$d(x_1, Tx_1) < \max_{i=1,...,k} g_i(M_{x_0 x_1})$$

$$= \max_{i=1,...,k} g_i(d(x_0, x_1), d(x_1, Tx_1), d(x_0, Tx_0), d(x_0, Tx_1), d(x_1, Tx_0))$$

$$= \max_{i=1,\ldots,k} g_i(d(x_0,x_1), d(x_1,Tx_1), d(x_0,x_1), d(x_0,Tx_1), 0)$$

$$\leq \max_{i=1,\ldots,k} g_i(d(x_0,x_1), d(x_1,Tx_1), d(x_0,x_1), d(x_0,x_1)+d(x_1,Tx_1), 0).$$

By proposition 1.11, we obtain $d(x_1, Tx_1) < d(x_0, x_1)$. Now let $x_2 = Tx_1$. Since $\alpha d(x_1, Tx_1) < d(x_1, x_2)$, by using the assumption we have

$$\psi[d(Tx_1, Tx_2)] \leq \phi(g_1(M_{x_1 x_2}), g_2(M_{x_1 x_2}), \ldots, g_k(M_{x_1 x_2})).$$

Then by lemma 1.12 we have

$$d(x_2, Tx_2) < \max_{i=1,\ldots,k} g_i(M_{x_1 x_2})$$

$$= \max_{i=1,\ldots,k} g_i(d(x_1,x_2), d(x_2,Tx_2), d(x_1,Tx_1), d(x_1,Tx_2), d(x_2,Tx_1))$$

$$= \max_{i=1,\ldots,k} g_i(d(x_1,x_2), d(x_2,Tx_2), d(x_1,x_2), d(x_1,Tx_2), 0)$$

$$\leq \max_{i=1,\ldots,k} g_i(d(x_1,x_2), d(x_2,Tx_2), d(x_1,x_2), d(x_1,x_2)+d(x_2,Tx_2), 0).$$

By proposition 1.11, we obtain $d(x_2, Tx_2) < d(x_1, x_2)$. Now by continuing this process, we obtain a sequence $\{x_n\}_{n\geq 1}$ in $X$ such that $x_{n+1} = Tx_n$ and $d(x_n, x_{n+1}) < d(x_{n-1}, x_n)$. So by lemma 1.13, there is $a$ such that $\lim_{n\to\infty} d(x_n, x_{n+1}) = a$. Hence

$$\lim_{n\to\infty} M_{x_n x_{n+1}} = (a, a, a, A, 0)$$

where $A \leq 2a$. Then

$$\psi(a) = \lim_{n\to\infty} \psi[d(x_{n+1}, x_{n+2})]$$

$$\leq \lim_{n\to\infty} \phi(g_1(M_{x_n x_{n+1}}), g_2(M_{x_n x_{n+1}}), \ldots, g_k(M_{x_n x_{n+1}}))$$

$$\leq \phi(g_1(a,a,a,2A,0), g_2(a,a,a,2A,0), \ldots, g_k(a,a,a,2A,0)).$$

Now by lemma 1.12 we obtain

$$a < \max_{i=1,\ldots,k} g_i(a,a,a,A,0) \leq \max_{i=1,\ldots,k} g_i(a,a,a,2a,0) \leq a$$

an then $a = 0$. We claim that $\{x_n\}_{n\geq 1}$ is a Cauchy sequence in $(X, d)$. Suppose that $\{x_n\}_{n\geq 1}$ is not a Cauchy sequence, which means that there is a constant $\varepsilon_0 > 0$ such that for each positive integer $k$, there are positive integers $m(k)$ and $n(k)$ with $m(k) > n(k) > k$ such that

$$d(x_{m(k)}, x_{n(k)}) \geq \varepsilon_0, \ d(x_{m(k)1}, x_{n(k)}) < \varepsilon_0.$$

From Lemma 1.14, we obtain,

$$\lim_{k\to\infty} d(x_{m(k)+1}, x_{n(k)+1}) = \varepsilon_0 \quad and \quad \lim_{k\to\infty} d(x_{m(k)+2}, x_{n(k)+2}) = \varepsilon_0.$$

We claim that for any $y \in X$, one of the flowing relations is held:
$$\alpha d(x_n, Tx_n) \leq d(x_n, y) \quad or \quad \alpha d(x_{n+1}, Tx_{n+1}) \leq d(x_{n+1}, y). \tag{1}$$
Otherwise, if $\alpha d(x_n, Tx_n) > d(x_n, y)$ and $\alpha d(x_{n+1}, Tx_{n+1}) > d(x_{n+1}, y)$, we have
$$d(x_n, x_{n+1}) \leq d(x_n, y) + d(x_{n+1}, y) < \alpha d(x_n, Tx_n) + \alpha d(x_{n+1}, Tx_{n+1})$$
$$= \alpha d(x_n, x_{n+1}) + \alpha d(x_{n+1}, x_{n+2}) \leq 2\alpha d(x_n, x_{n+1}) \leq d(x_n, x_{n+1})$$
which is a contradiction. Now by using the assumption and relation 1, for each $n \geq 1$ one of the following cases holds:
(i) There exists an infinite subset $I \subset \mathbb{N}$ such that
$$\psi[d(x_{m(k)+1}, x_{n(k)+1})]$$
$$\leq \phi(g_1(M_{x_{m(k)} x_{n(k)}}), g_2(M_{x_{m(k)} x_{n(k)}}), ..., g_k(M_{x_{m(k)} x_{n(k)}})).$$
(ii)There exists an infinite subset $J \subset \mathbb{N}$ such that
$$\psi[d(x_{m(k)+2}, x_{n(k)+1})]$$
$$\leq \phi(g_1(M_{x_{m(k)+1} x_{n(k)}}), g_2(M_{x_{m(k)+1} x_{n(k)}}), ..., g_k(M_{x_{m(k)+1} x_{n(k)}})).$$
Since
$$M_{x_{m(k)} x_{n(k)}}$$
$$= (d(x_{m(k)}, x_{n(k)}), d(x_{n(k)}, Tx_{n(k)}), d(x_{m(k)}, Tx_{m(k)}), d(x_{m(k)}, Tx_{n(k)}), d(x_{n(k)}, Tx_{m(k)}))$$
$$= (d(x_{m(k)}, x_{n(k)}), d(x_{n(k)}, x_{n(k)+1}), d(x_{m(k)}, x_{m(k)+1}), d(x_{m(k)}, x_{n(k)+1}), d(x_{n(k)}, x_{m(k)+1}))$$
$$\leq (d(x_{m(k)}, x_{n(k)}), d(x_{n(k)}, x_{n(k)+1}), d(x_{m(k)}, x_{m(k)+1}),$$
$$d(x_{m(k)}, x_{n(k)}) + d(x_{n(k)}, x_{n(k)+1}), d(x_{n(k)}, x_{m(k)}) + d(x_{m(k)}, x_{m(k)+1})),$$
we have $\lim_{n \to \infty} M_{x_{m(k)} x_{n(k)}} = (\varepsilon_0, 0, 0, A, B)$ where $A \leq \varepsilon_0$ and $B \leq \varepsilon_0$. Then in case (i), we obtain
$$\psi(\varepsilon_0)$$
$$\leq \phi(g_1(\varepsilon_0, 0, 0, A, B), g_2(\varepsilon_0, 0, 0, A, B), ..., g_k(\varepsilon_0, 0, 0, A, B))$$
and then by lemma 1.12 we have
$$\varepsilon_0 < \max_{i=1,...,k} g_i(\varepsilon_0, 0, 0, A, B) \leq \max_{i=1,...,k} g_i(\varepsilon_0, 0, 0, \varepsilon_0, \varepsilon_0) \leq \varepsilon_0,$$
which is a contradiction.
In case (ii), similar to cas(i), we obtain
$$\varepsilon_0 < \varepsilon_0,$$
which is a contradiction. This proves our claim that $\{x_n\}_{n \geq 1}$ is a Cauchy sequence in $(X, d)$. Let $\lim_{n \to \infty} x_n = x$. By relation 1, for each $n \geq 1$ and $y \in X$, either
a) $\psi[d(Tx_n, Ty)] \leq \phi(g_1(M_{x_n x}), g_2(M_{x_n x}), ..., g_k(M_{x_n x}))$

or

b) $\psi[d(Tx_{n+1}, Ty)] \leq \phi(g_1(M_{x_{n+1}x}), g_2(M_{x_{n+1}x}), ..., g_k(M_{x_{n+1}x}))$

In case (a), by using of lemma 1.12 we obtain

$$d(x, Tx) \leq d(x, Tx_n) + d(Tx_n, Tx) < d(x, Tx_n) + \max_{i=1,...,k} g_i(M_{x_n x})$$

$$= d(x, Tx_n) + \max_{i=1,...,k} g_i(d(x_n, x), d(x_n, Tx_n), d(x, Tx), d(x, Tx_n), d(x_n, Tx)).$$

Hence

$$d(x, Tx) \leq \max_{i=1,...,k} g_i(0, 0, d(x, Tx), 0, d(x, Tx)).$$

Now by using Proposition 1.11, we have $d(x, Tx) = 0$ and so $x = Tx$.

In case (b), by using lemma 1.12, we obtain

$$d(x, Tx) \leq d(x, Tx_{n+1}) + d(Tx_{n+1}, Tx) < d(x, Tx_n) + \max_{i=1,...,k} g_i(M_{x_{n+1}x})$$

$$\leq d(x, Tx_{n+1}) + \max_{i=1,...,k} g_i(d(x_{n+1}, x), d(x_{n+1}, Tx_{n+1}), d(x, Tx), d(x, Tx_{n+1}), d(x_{n+1}, Tx)).$$

Hence

$$d(x, Tx) \leq g(0, 0, d(x, Tx), 0, d(x, Tx)),$$

and then by using Proposition 1.11, we have $d(x, Tx) = 0$. So $x = Tx$. We claim that this fixed point is unique. Suppose that there are two distinct points $a, b \in X$ such that $Ta = a$ and $Tb = b$. Since $d(a, b) > 0 = \alpha d(a, Ta)$, we have the contradiction

$$0 < \psi[d(a, b)] = \psi[d(Ta, Tb)]$$

$$\leq \phi(g_1(M_{ab}), g_2(M_{ab}), ..., g_k(M_{ab})).$$

Now by lemma 1.12, we obtain

$$d(a, b) < \max_{i=1,...,k} g_i(d(a, b), d(a, Ta), d(b, Tb), d(a, Tb), d(b, Ta))$$

$$= \max_{i=1,...,k} g_i(d(a, b), 0, 0, d(a, b), d(b, a)) \leq d(a, b).$$

So $d(a, b) = 0$. $\square$

**Corollary 2.2.** *Let $(X, d)$ be a complete metric space and $T : X \to X$ be a mapping satisfying*

$$\psi(d(Tx, Ty)) \leq h(\psi(M(x, y)))$$

*for each $x, y \in X$, where*

$$M(x, y) = \max\{d(x, y), d(x, Tx), d(y, Ty), d(x, Ty), d(y, Tx)\},$$

*$\psi \in \Psi$ and $h : \mathbb{R}^+ \to \mathbb{R}^+$ is a continuous function such that $h(t) < t$ for all $t > 0$. Then $T$ has a unique fixed point $x \in X$.*

*Proof.* Let $g_1(t_1,t_2,t_3,t_4,t_5) = \max\{t_1,t_2,t_3,t_4,t_5\}$ and define $\phi$ by $\phi(t) = h(\psi(t))$. It is easy to see that $\phi \in \Phi$ and for every $t > 0$, $\phi(t) < \psi(t)$. Now by using Theorem 2.1, $T$ has a fixed point. □

**Remark 2.1.** *By lemma 1.5, we see that theorems 1.2, 1.3 and 1.4 are special cases of theorem 2.1.*

**Remark 2.2.** *Theorem 1.7 is a special case of theorem 2.1.*

*Proof.* Let $g_1 = g_2(t_1,t_2,t_3,t_4,t_5) = t_1$ and define $\phi$ by $\phi(t_1,t_2) = \phi(t_1) - h(t_2)$. Now by using Theorem 2.1, $T$ has a fixed point. □

**Remark 2.3.** *Theorem 1.8 is a special case of theorem 2.1.*

*Proof.* Let $g_1(t_1,t_2,t_3,t_4,t_5) = \max\{t_1,t_2,t_3,\frac{1}{2}(t_4+t_5)\}$ and define $\phi(t) = \psi(t) - h(t)$. Now by using Theorem 2.1, $T$ has a fixed point. □

**Remark 2.4.** *Theorem 1.9 is a special case of theorem 2.1.*

*Proof.* Let $g_1(t_1,t_2,t_3,t_4,t_5) = \max\{t_1,t_2,t_3,\frac{1}{2}(t_4+t_5)\}$, $g_2(t_1,t_2,t_3,t_4,t_5) = \max\{t_1,t_2\}$ and define $\phi(t_1,t_2) = \psi(t_1) - h(t_2)$. Now by using theorem 2.1, $T$ has a fixed point. □

**Remark 2.5.** *Let $g_1(t_1,t_2,t_3,t_4,t_5) = t_1$, $g_2(t_1,t_2,t_3,t_4,t_5) = t_2\frac{1+t_3}{1+t_1}$ and define $\phi(t_1,t_2) = a\psi(t_1) + b\psi(t_2)$. Then we obtain theorem 1.10 of theorem 2.1.*